# A COMPLETE DESCRIPTION OF THE ASYMPTOTIC BEHAVIOR AT INFINITY OF POSITIVE RADIAL SOLUTIONS TO $\Delta^2 u = u^\alpha$ IN $\mathbf{R}^n$

QUỐC ANH NGÔ, VAN HOANG NGUYEN, AND QUOC HUNG PHAN

Abstract. We consider the biharmonic equation $\Delta^2 u = u^\alpha$ in $\mathbf{R}^n$ with $n \geq 1$. It was proved that this equation has a positive classical solution if, and only if, either $\alpha \leq 1$ with $n \geq 1$ or $\alpha \geq (n+4)/(n-4)$ with $n \geq 5$. The asymptotic behavior at infinity of all positive radial solutions was known in the case $\alpha \geq (n+4)/(n-4)$ and $n \geq 5$. In this paper, we classify the asymptotic behavior at infinity of all positive radial solutions in the remaining case $\alpha \leq 1$ with $n \geq 1$; hence obtaining a complete picture of the asymptotic behavior at infinity of positive radial solutions. Since the underlying equation is higher-order, we propose a new approach that relies on a representation formula and asymptotic analysis arguments. We believe that the approach introduced here can be conveniently applied to study other problems with higher-order operators.

## 1. Introduction

The aim of the present paper is to study positive classical solutions to the following biharmonic equation

$$\Delta^2 u = u^\alpha \tag{1.1}$$

in the whole Euclidean space $\mathbf{R}^n$ with $n \geq 1$ and $\alpha \in \mathbf{R}$. Such the equation has already captured so much attention in the last two decades, that will be described later. From now on, we shall call a classical solution on the whole space $\mathbf{R}^n$ an entire solution.

The motivation of working on (1.1) comes from a rapidly increasing number of papers on higher order elliptic equations in $\mathbf{R}^n$ in recent years. The biharmonic equation (1.1) is a higher-order analog of the Lane-Emden equation

$$-\Delta u = u^\alpha \tag{1.2}$$

in $\mathbf{R}^n$, which has already been in the core of many researches in the last few decades. Concerning (1.2), there is a threshold $p_{\mathrm{S}}^2$, known as the critical Sobolev exponent, which is given as follows

$$p_{\mathrm{S}}^2 = \begin{cases} (n+2)/(n-2) & \text{if } n \geq 3, \\ +\infty & \text{if } n \leq 2. \end{cases}$$

A fundamental existence result tells us that an entire positive solution to (1.2) exists if and only if $n \geq 3$ and $\alpha \geq p_{\mathrm{S}}^2$; see [GS81, JL73, AS11] and references therein.

Back to (1.1), the first question is to know under what conditions on $\alpha$ an entire solution does actually exist. As far as we know, a complete answer has been established. To be precise, it was proved by Lin [Lin98] that there is no entire positive solution to

---

(1.1) whenever $1 < \alpha < p_S^4$, where $p_S^4$ is the other critical Sobolev exponent defined as follows

$$p_S^4 = \begin{cases} (n+4)/(n-4) & \text{if } n \geq 5, \\ +\infty & \text{if } n \leq 4. \end{cases}$$

In the critical and super-critical cases, namely $\alpha = p_S^4$ and $\alpha > p_S^4$ respectively, the existence of entire positive solutions to (1.1) was respectively shown by Lin [**Lin98**] and Gazzola and Grunau [**GG06**]. In the remaining case $\alpha \in (-\infty, 1]$, not long ago, it was proved in [**NNPY20**, Proposition 4.5] that (1.1) has at least a positive radial solution.

The motivation of writing this paper traces back to the works of Lin and of Gazzola and Grunau mentioned above. In the critical case, a beautiful classification result in [**Lin98**] indicates that any entire positive solution to (1.1) are of the form

$$u(x) = \Big(\frac{2\lambda}{1+\lambda^2|x-x_0|^2}\Big)^{\frac{n-4}{2}}.$$

In the super-critical case, among other things, the authors in [**GG06**] showed that all entire positive radial solutions to (1.1) for $\alpha > p_S^4$ obey the following asymptotic behavior

$$\lim_{r \nearrow +\infty} r^{-\frac{4}{1-\alpha}} u(r) = L,$$

where

$$L = \Big[\big(m(m+2)(n-2-m)(n-4-m)\Big]^{m/4}$$

with $m := 4/(\alpha-1)$. For further understanding on the asymptotic behavior of radial solutions to (1.1), we refer to [**GG06**, **FGK09**, **Kar09**, **Win10**]. Therefore, inspired by the results obtained in [**Lin98**] and [**GG06**] for $\alpha \geq p_S^4$, in this paper, as a counterpart, we focus our attention on the asymptotic behavior of positive radial solutions to (1.1) in the range $\alpha \leq 1$ for any dimension $n \geq 1$. To be more precise, our primary aim is to classify the growth and the asymptotic behavior at infinity of any positive radial solutions to (1.1), thus completing the picture of the asymptotic behavior at infinity of positive radial solutions of problem (1.1) in $\mathbf{R}^n$. However, we emphasize that unlike the critical and super-critical range, the asymptotic behavior at infinity for $\alpha \leq 1$ contains many generic cases; see Table 1 below. In order to let the finding easily accessible, there will be no limit for the dimension $n$, namely the result obtained in this paper is valid for any $n \geq 1$. Furthermore, there will also be no additional assumption on $\alpha$ and on radial solutions.

Before closing this section, we would like to mention that seeking for the asymptotic behavior at infinity for positive solutions to partial differential equations is a classic question. In the analogous second-order problem (1.2), the asymptotic behavior at infinity for positive solutions to (1.2) was completely classified by Ni in [**Ni82**] for the critical case and by Wang in [**Wan93**] for the super-critical case. In the super-critical case, Wang found that any radial solution to (1.2) obeys the following asymptotic behavior

$$\lim_{r \nearrow +\infty} r^{-\frac{2}{1-\alpha}} u(r) = \Big[\frac{2}{\alpha-1}\Big(n-2-\frac{2}{\alpha-1}\Big)\Big]^{-\frac{1}{1-\alpha}}.$$

A counterpart of (1.2) is the following equation

$$\Delta u = u^\alpha \tag{1.3}$$

in $\mathbf{R}^n$ and in this scenario, it is easily deduced that a positive solution to (1.3) exists provided $\alpha \leq 1$; also see [**NNPY20**] for further existence and non-existence results. Under the dimensional restriction $n \geq 3$ and as far as we know, the asymptotic behavior of positive radial solutions to equation (1.3) was studied by Yang and Guo in [**YG05**] for the case $\alpha \in (0, 1)$ and by Guo, Guo, and Li in [**GGL06**] for the case $\alpha < 0$. Combining the results obtained in [YG05] and [GGL06], it is now known that any radial solution to (1.3) for $n \geq 3$ and $\alpha < 1$ obeys the following asymptotic behavior at infinity

$$\lim_{r \nearrow +\infty} r^{-\frac{2}{1-\alpha}} u(r) = \left[\frac{2}{1-\alpha}\left(n-2+\frac{2}{1-\alpha}\right)\right]^{-\frac{1}{1-\alpha}}.$$

Still for the asymptotic behavior of radial solutions to (1.3) when $n \geq 3$, the only case left is $\alpha = 1$. It is expected that such a result should be well-known, however, we are unable to find a reference for it. Therefore, we manage to calculate the asymptotic behavior in this remaining case and put it in Appendix A; see Proposition A.1. In the same spirit, the following counterpart of (1.1)

$$\Delta^2 u = -u^\alpha \tag{1.4}$$

in $\mathbf{R}^n$ has attracted much attention starting from the preliminary version of a paper by Choi and Xu [CX09] and a paper by McKenna and Reichel [MR03]. In this case, it has been proved that (1.4) has at least a positive solution if and only if $\alpha < -1$ and $n \geq 3$; see [LY16, NNPY20]. Taking this restriction into account, the asymptotic behavior at infinity of positive radial solutions to (1.4) has been widely investigated, for instance in [GW08, DFG10, Gue12, Wan14, DN17].

As we are moving toward the end of this section, we would like to draw some direction for future research. It is worth emphasizing that the description of the asymptotic behavior of radial solutions to (1.1) is the first step in understanding the picture of the set of solutions to the equations. Unlike the counterpart (1.4), equation (1.1) in the range $\alpha \leq 1$ always admits at least one 'radial' solution; see [NNPY20, Proposition 4.5]. Therefore, it is easy to produce a new non-radial solution to (1.1) in higher dimensions. Let us discuss this in detail. Suppose that $u(x)$ is a positive radial solution to (1.1) in $\mathbf{R}^n$, namely

$$\Delta_x^2 u(x) = u(x)^\alpha$$

in $\mathbf{R}^n$. Then the new function

$$v(x,y) := u(x)$$

also solves (1.1) in $\mathbf{R}^{n+m}$ for any $m \geq 1$, namely

$$\Delta_{(x,y)}^2 v(x,y) = v(x,y)^\alpha$$

in $\mathbf{R}^{n+m}$. It is clear that $v(x,y)$ is no longer radially symmetric in $\mathbf{R}^{n+m}$. As can be easily seen from Table 1, radial solutions to (1.1) grows at least linearly at infinity. Although we are not sure if this is still true for non-radial solutions, but if this is the case, the growth of $u$ at infinity in $\mathbf{R}^n$ is the maximal growth of $v$ at infinity in $\mathbf{R}^{n+m}$. Of course, it is not clear if there are other non-radial solutions to (1.1) rather than the ones described above. As far as we know, even for the lower order case, existence of non-symmetric solutions to $\Delta u = u^\alpha$ in $\mathbf{R}^n$ with $0 < \alpha < 1$ remains open. By non-symmetry we mean that the solution is not symmetric with respect to any line passing through the origin in the ambient space. Using this convention, it is easy to realize that any radial function $u$ in $\mathbf{R}^n$ is still symmetric with respect to the $u$-axis in $\mathbf{R}^{n+1}$. Nevertheless, toward a complete picture of solutions to (1.1), it is no doubt that the asymptotic behavior obtained in this paper will play some role. Knowing that the fourth order Hardy–Hénon equation,

$$\Delta^2 u = |x|^\sigma u^\alpha \quad \text{in } \mathbf{R}^n \tag{1.5}$$

with $\sigma \in \mathbf{R}$ has already been captured attention, we hope that our finding in this work corresponding to $\sigma = 0$ can be generalized to $\sigma \neq 0$. It is worth noting that the two cases $\alpha > 1$ and $\alpha \leq 1$ are quite different. When $\alpha > 1$ there are several results concerning the existence and non-existence of solutions to (1.5), see [DPQ22, NY22] for the case $n \geq 5$ and [NNT24, Ngo24] for the case $2 \leq n \leq 4$. However, in the case $\alpha \leq 1$, the picture of existence and non-existence of solutions to (1.5) is not clear. As far as we know, the work [DQ20] by Dai and Qin is the first one considering $0 < \alpha \leq 1$. Further handling $\sigma \neq 0$ and $\alpha \leq 1$ is our work in future.

In the next section, we state our main results whose proofs are put in Section 3. But the precise structure of this paper is organized as follows.

Contents

Finally, as a remark before closing this section, we should mention that there are arguments and results more or less known to experts in this field. In fact, we believe that some of the findings can be argued differently using certain known results. However, we aim to argue as much elementary and simple as possible, for the reader's convenience, while trying to maintain the paper in a reasonable length.

## 2. Statement of main results

For convenience, by the notation $u \sim f$ we mean that

$$\lim_{r \nearrow +\infty} \frac{u(r)}{f(r)} = 1.$$

Sometime, to emphasize that the above limit is taking as $r \searrow 0$, we shall write $\sim_{r \searrow 0}$. Let us denote the following universal constant

$$C(n, \alpha) = \Big(n + \frac{4\alpha}{1-\alpha}\Big)\Big(n + \frac{4\alpha}{1-\alpha} + 2\Big)\Big(\frac{4\alpha}{1-\alpha} + 2\Big)\Big(\frac{4\alpha}{1-\alpha} + 4\Big).$$

In $\mathbf{R}^n$, it is well-known that the Laplace and biLaplace operators acting on a radial function $u$ can be expressed as follows

$$\Delta u(r) = \frac{1}{r^{n-1}}\big(r^{n-1} u'(r)\big)' = u''(r) + \frac{n-1}{r} u'(r) \tag{2.1}$$

and

$$\Delta^2 u(r) = u^{(4)}(r) + \frac{2(n-1)}{r} u^{(3)}(r) + \frac{(n-1)(n-3)}{r^2} u''(r) - \frac{(n-1)(n-3)}{r^3} u'(r) \tag{2.2}$$

for $r > 0$. It is well-known that

$$\Delta^2 \left( r^{\frac{4}{1-\alpha}} \right) = C(n,\alpha) r^{\frac{4\alpha}{1-\alpha}}.$$

We are now in position to state our main results. To be more precise, we shall provide an asymptotic expansion for any radial solution to (1.1) near infinity. It is worth emphasizing that our results require no condition on $n \geq 1$ and on $\alpha$ except $\alpha \leq 1$, which is natural based on the discussion described in Introduction. However, since the formulation of the results are rather long and the technique used is different when $n$ varies, we intend to split our results into three theorems according to either the dimension $n \geq 3$, $n = 2$, or $n = 1$.

First, for the case $n \geq 3$, our result reads as follows.

**Theorem 1.** *Assume $n \geq 3$ and let $u$ be a positive radial solution to the problem (1.1) in $\mathbf{R}^n$. Then we have the following claims:*

*(a) If $\alpha = 1$, then*

$$u \sim \Big[ u(0) + \Delta u(0) \Big] \Gamma\Big(\frac{n}{2}\Big) 2^{\frac{n-5}{2}} \pi^{-1/2} r^{-\frac{n-1}{2}} e^r. \tag{2.3}$$

*(b) If $\alpha \in (-1,1)$, then*

$$u \sim C(n,\alpha)^{-\frac{1}{1-\alpha}} r^{\frac{4}{1-\alpha}}. \tag{2.4}$$

*(c) If $\alpha = -1$, then*

$$u \sim \Big( n(n-2) \Big)^{-1/2} r^2 (\log r)^{1/2}. \tag{2.5}$$

*(d) If $\alpha < -1$, then*

$$u \sim \frac{1}{2n} \Big( \Delta u(0) + \frac{1}{n-2} \int_0^{+\infty} t u^\alpha(t) dt \Big) r^2. \tag{2.6}$$

When the dimension $n = 2$, we obtain the following asymptotic expansion near infinity.

**Theorem 2.** *Let $u$ be a positive radial solution to the problem (1.1) in $\mathbf{R}^2$, then we have the following claim:*

*(a) If $\alpha = 1$, then*

$$u \sim \Big[ u(0) + \Delta u(0) \Big] 2^{-3/2} \pi^{-1/2} r^{-1/2} e^r.$$

*(b) If $\alpha \in (-1,1)$, then*

$$u \sim C(2,\alpha)^{-\frac{1}{1-\alpha}} r^{\frac{4}{1-\alpha}}.$$

*(c) If $\alpha = -1$, then*

$$u \sim 2^{-1/2} r^2 \log r (\log\log r)^{1/2}. \tag{2.7}$$

*(d) If $\alpha < -1$, then*

$$u \sim \Big( \frac{1}{4} \int_0^{+\infty} t u^\alpha(t) dt \Big) r^2 \log r. \tag{2.8}$$

We note that the asymptotic behavior in Theorem 2(a) is exactly the same as that of (2.3) with $n$ replaced by 2 and $\alpha$ replaced by 1. Finally, our result for dimension $n = 1$ is as follows.

**Theorem 3.** *Let $u$ be a positive radial solution to the problem (1.1) in $\mathbf{R}$, then we have the following claim:*

(a) *If $\alpha = 1$, then*
$$u \sim \frac{1}{4}\Big[u(0) + u''(0)\Big]e^r.$$
(b) *If $\alpha \in (-1/3, 1)$, then*
$$u \sim C(1,\alpha)^{-\frac{1}{1-\alpha}} r^{\frac{4}{1-\alpha}}.$$
(c) *If $\alpha = -1/3$, then*
$$u \sim \Big(\frac{2}{9}\Big)^{3/4} r^3 (\log r)^{3/4}. \tag{2.9}$$
(d) *If $\alpha < -1/3$, then*
$$u \sim \frac{1}{6}\Big(\int_0^{+\infty} u^\alpha(t)dt\Big) r^3. \tag{2.10}$$

Again, we note that the asymptotic behavior in Theorem 3(a) is exactly the same as that of (2.3) with both $n$ and $\alpha$ replaced by 1.

As can be easily notified, the asymptotic behavior of radial solutions to (1.1) in the case $-1 \leq \alpha < 1$ with $n \geq 2$ and in the case $-1/3 \leq \alpha < 1$ with $n = 1$ does not depend on the solution.

The proof of our main results follows from a general procedure, which is based on a suitable combination of a priori bounds from below and above, integral estimates and asymptotic analysis arguments.

The following table summarizes the asymptotic behavior at infinity of entire positive radial solutions to (1.1) without the precise limit and gives the sketch of proof of Theorems 1-3. (By comparing with the case $\alpha > 1$, one easily notices that the case $\alpha \leq 1$ is much more completed; see [**Lin98**] and [**GG06**].)

<table>
<tr><th>$n$</th><th>$\alpha < -1$</th><th>$\alpha = -1$</th><th>$-1 < \alpha < -\frac{1}{3}$</th><th>$\alpha = -\frac{1}{3}$</th><th>$-\frac{1}{3} < \alpha < 1$</th><th>$\alpha = 1$</th></tr>
<tr><td>1</td><td colspan="3">$r^3$<br>Prop. 3.8</td><td>$r^3(\log r)^{3/4}$<br>Prop. 3.7</td><td>$r^{\frac{4}{1-\alpha}}$<br>Prop. 3.2</td><td rowspan="3">$r^{-\frac{n-1}{2}} e^r$<br>Prop. 3.1</td></tr>
<tr><td>2</td><td>$r^2 \log r$<br>Prop. 3.6</td><td>$r^2 \log r \sqrt{\log \log r}$<br>Prop. 3.5</td><td colspan="3" rowspan="2">$r^{\frac{4}{1-\alpha}}$<br>Prop. 3.2</td></tr>
<tr><td>$\geq 3$</td><td>$r^2$<br>Prop. 3.4</td><td>$r^2 \sqrt{\log r}$<br>Prop. 3.3</td></tr>
</table>

Table 1: Asymptotic behavior at infinity of entire positive radial solutions to $\Delta^2 u = u^\alpha$ in $\mathbf{R}^n$ with $n \geq 1$ and $\alpha \leq 1$.

From Table 1 above, all radial solutions to (1.1) are *sky state*, a terminology introduced in [**GGL06**]. Moreover, radial solutions tend to grow stronger when $\alpha$ increases toward 1 until $\alpha$ reaches 1 where a non-existence result occurs. Combining with the asymptotic behavior obtained in [**Lin98**] and [**GG06**], after passing through the non-existence regime, radial solutions tend to decay and they decay more rapidly as $\alpha$ increases toward infinity.

Before closing this section, we note that in the case $\alpha < -1$ with $n \geq 2$, Theorems 1(d) and 2(d) were partially proved by Kusano, Naito, and Swanson in [KNS87, KNS88]. More precisely, it was showed in [KNS87, Theorem 2] that the quotient $u(r)/(r^2 \log r)$ has a limit as $r \nearrow +\infty$ when $n = 2$ and $\alpha < -1$. When $n \geq 3$ and $\alpha < -1$ it was proved in [KNS88, Theorem 2] that $u(r)/r^2$ has a limit as $r \nearrow +\infty$. However, the exact value of these limits were not computed in these mentioned works. Our finding clearly provides the exact value for these limits.

## 3. Proofs of the main results

This section is devoted to the proof of our main result. For the sake of clarity, we divide this section into several parts. We spend Subsection 3.1 to collect some auxiliary results while proofs for Theorems 1–3 are put in Subsection 3.4–3.6.

### 3.1. Some useful estimates.

In this subsection, we collect some basic results which are used many times in our arguments.

**Lemma 3.1.** *Let $v$ be a radial $C^2$-function in $\mathbf{R}^n$ with $n \geq 1$. Then we have the following representation*

$$v(r) = v(r_0) + \int_{r_0}^r s^{1-n}\Big(\int_0^s t^{n-1}\Delta v(t)dt\Big)ds$$

*and*

$$v(r) = v(r_0) + \int_0^{r_0} t^{n-1}\Delta v(t)dt \int_{r_0}^r s^{1-n}ds + \int_{r_0}^r s^{1-n}\Big(\int_{r_0}^s t^{n-1}\Delta v(t)dt\Big)ds$$

*for any fixed $r_0 \geq 0$.*

*Proof.* This is elementary. For the first identity, suppose $n = 1$. Since $v$ is radial, there holds $v'(0) = 0$. Hence, the identity follows from $v'' = \Delta v$ via integration by parts. When $n \geq 2$, we integrate both sides of

$$(r^{n-1}v'(r))' = r^{n-1}\Delta v(r)$$

over $[r_0, r]$ to get the desired identity. For the second identity, for $n \geq 2$, this identity comes from the first identity by splitting the domain of integration. When $n = 1$, the identity

$$v(r) = v(r_0) + (r - r_0)\int_0^{r_0} \Delta v(t)dt + \int_{r_0}^r \int_{r_0}^s \Delta v(t)dtds$$

follows from integration by parts with a note that $v'(r_0) = \int_0^{r_0} v''(t)dt$. This completes the proof. ☐

Our next useful estimate is about the growth of radial solutions to (1.1).

**Lemma 3.2.** *Let $u$ be a positive radial solution to the problem (1.1) in $\mathbf{R}^n$ with $\alpha < 1$. Then, there holds*

$$\lim_{r \nearrow +\infty} \Delta u(r) > 0.$$

*Consequently, there exist two positive numbers $C_0$ and $r_0$ such that*

$$u(r) \geq C_0 r^2$$

*for all* $r \geq r_0$.

*Proof.* First we observe that

$$\left(r^{n-1}(\Delta u)'(r)\right)' = r^{n-1}\Delta^2 u > 0,$$

which implies that the function $\Delta u$ is increasing. Hence, there exists

$$\gamma := \lim_{r \nearrow +\infty} \Delta u(r),$$

which could be infinity. Suppose that $\gamma \leq 0$. Then $\Delta u(r) \leq 0$ for any $r \geq 0$. Since

$$\Delta u(r) = r^{1-n}(r^{n-1}u'(r))'$$

we deduce that $u(r)$ is non-increasing. From this it follows that $u(r) \leq u(0)$ for any $r \geq 0$. Depending on the value of $\alpha$, we consider the following two cases:

**Case 1**. Suppose that $\alpha \in [0,1)$. We deduce from the rescaled test-function argument that

$$\int_{B_R} u^\alpha dx \leq CR^{-4} \int_{B_{2R}} u dx \leq Cu^{1-\alpha}(0)R^{-4} \int_{B_{2R}} u^\alpha dx.$$

Repeating this argument, we obtain, for any $m \geq 1$, the following estimate

$$\int_{B_R} u^\alpha dx \leq \left(Cu^{1-\alpha}(0)\right)^m R^{-4m} \int_{B_{2^m R}} u^\alpha dx.$$

Now choosing any integer $m > n/4$ and using the fact that $u(r) \leq u(0)$, we have

$$\int_{B_R} u^\alpha dx \leq CR^{-4m}(2^m R)^n \leq CR^{n-4m}.$$

From this we let $R \nearrow +\infty$ to get $\int_{\mathbf{R}^n} u^\alpha dx = 0$, which is impossible.

**Case 2**. Suppose that $\alpha < 0$. In this case, we can apply Lemma 3.1 to get

$$\begin{aligned}
\Delta u(r) &= \Delta u(0) + \int_0^r s^{1-n} \int_0^s t^{n-1} u^\alpha(t) dt ds \\
&\geq \Delta u(0) + u^\alpha(0) \int_0^r s^{1-n} \int_R^s t^{n-1} dt ds \\
&\geq \Delta u(0) + \frac{u^\alpha(0)}{2n} r^2.
\end{aligned}$$

Hence, $\Delta u(r)$ is positive for $r$ large, which contradicts $\gamma \leq 0$.

Combining two cases above, we deduce $\gamma > 0$. From this one can easily obtain the existence of $C_0$ and $r_0$ as claimed. The proof is complete. □

Our last useful estimate is a consequence of the l'Hôpital rule, whose proof is omitted.

**Lemma 3.3**. *Let* $v \in C^2[t_0, +\infty)$ *for some* $t_0 > 0$. *Assume* $\sigma > -2$ *and*

$$\lim_{r \nearrow +\infty} \frac{\Delta v(r)}{r^\sigma} = L \in \mathbf{R}.$$

*Then there holds*

$$\lim_{r \nearrow +\infty} \frac{v(r)}{r^{\sigma+2}} = \frac{L}{(2+\sigma)(n+\sigma)}.$$

3.2. **The special case:** $\alpha = 1$ **with** $n \geq 1$. This subsection is devoted to a proof for part of Theorems 1–3 indicated in Table 1. We start with the special case $\alpha = 1$. This case is special because, recall (2.2), the equation $\Delta^2 u(r) = u(r)$ in $\mathbf{R}^n$ can be written in terms of the following ODE

$$u^{(4)}(r) + \frac{2(n-1)}{r} u^{(3)}(r) + \frac{(n-1)(n-3)}{r} u^{(2)}(r) - \frac{(n-1)(n-3)}{r} u'(r) = u(r) \qquad (3.1)$$

for $r > 0$ together with the following initial conditions

$$u'(0) = u'''(0) = 0.$$

**Proposition 3.1.** *Let $u$ be a positive solution to the ODE (3.1) with $n \geq 1$, then $u$ has the following asymptotic behavior*

$$u \sim \Big[u(0) + nu''(0)\Big] \Gamma\Big(\frac{n}{2}\Big) 2^{\frac{n-5}{2}} \pi^{-1/2} r^{-\frac{n-1}{2}} e^r$$

*near infinity.*

*Proof.* The argument below, relying heavily on the standard theory of ODE, is more or less known to experts in this field. However, we present a detailed argument for completeness. For clarity, the two cases $n = 1$ and $n \geq 2$ are considered separately.

**The case** $n = 1$. In this case, our equation (3.1) simply becomes the ODE

$$u^{(4)}(r) = u(r).$$

From this we choose four independent solutions $e^r$, $e^{-r}$, $\sin r$, and $\cos r$ to form the general solution to the ODE which is of the form

$$u(r) = C_1 e^r + C_2 e^{-r} + C_3 \sin r + C_4 \cos r$$

for suitable constants $C_i$ with $1 \leq i \leq 4$. These constants $C_i$ can be expressed in terms of $u(0)$ and $u^{(i)}(0)$ with $1 \leq i \leq 3$. To be more precise, there holds

$$C_1 = \frac{1}{4} \sum_{0 \leq i \leq 3} u^{(i)}(0),$$

$$C_2 = \frac{1}{4} \sum_{0 \leq i \leq 3} (-1)^i u^{(i)}(0),$$

$$C_3 = \frac{1}{2} [u'(0) - u^{(3)}(0)],$$

and

$$C_4 = \frac{1}{2} [u(0) - u^{(2)}(0)].$$

Note that, if $u$ is a positive solution of the form above, then it is necessary that $C_1 > 0$. In this case we obtain

$$u \sim C_1 e^r.$$

Keep in mind that $u'(0) = u^{(3)}(0) = 0$. Hence, we have just shown that

$$u \sim \frac{1}{4} \Big[u(0) + u''(0)\Big] e^r,$$

which is the desired behavior because $\Gamma(1/2) = \sqrt{\pi}$. A closer look at $C_i$ with $2 \leq i \leq 4$ reveals $C_3 = 0$ and therefore the general solution $u(r)$ to the ODE actually depends only on two parameters $u(0)$ and $u''(0)$, which can be written as follows

$$u(r) = \frac{u(0) + u''(0)}{4} \Big[e^r + e^{-r}\Big] + \frac{u(0) - u''(0)}{2} \cos r.$$

**The case** $n \geq 2$. As in the case $n = 1$, we need to find four independent solutions to the equation (3.1). With a little help from Maple applied to the formula for $\Delta^2 u(r)$ mentioned at the beginning of Section 2, we obtain the following four independent radial solutions

$$\begin{aligned} u_1(r) &= r^{-\frac{n}{2}+1} I_{\frac{n}{2}-1}(r), \\ u_2(r) &= r^{-\frac{n}{2}+1} K_{\frac{n}{2}-1}(r), \\ u_3(r) &= r^{-\frac{n}{2}+1} J_{\frac{n}{2}-1}(r), \end{aligned}$$

and

$$u_4(r) = r^{-\frac{n}{2}+1} Y_{\frac{n}{2}-1}(r),$$

where $I$ and $K$ are the modified Bessel functions of the first and second kinds respectively while $J$ and $Y$ are the Bessel functions of the first and second kinds respectively; see [AS64, Chapter 9]. Hence the general solution to the ODE is of the form

$$u(r) = \sum_{1 \leq i \leq 4} C_i u_i(r),$$

where $C_i$, $1 \leq i \leq 4$, are constants. Recall the asymptotic behavior for the Bessel and modified Bessel functions at infinity

$$\begin{aligned} I_{\frac{n}{2}-1}(r) &\sim \frac{e^r}{\sqrt{2\pi r}}, \\ K_{\frac{n}{2}-1}(r) &\sim \frac{\sqrt{\pi} e^{-r}}{\sqrt{2r}}, \\ J_{\frac{n}{2}-1}(r) &\sim \sqrt{\frac{2}{\pi r}} \cos\Big(r - \Big(\frac{n}{2} - 1\Big)\frac{\pi}{2} - \frac{\pi}{4}\Big), \end{aligned}$$

and

$$Y_{\frac{n}{2}-1}(r) \sim \sqrt{\frac{2}{\pi r}} \sin\Big(r - \Big(\frac{n}{2} - 1\Big)\frac{\pi}{2} - \frac{\pi}{4}\Big);$$

see, e.g., [AS64, Sections 9.2 and 9.7]. Again, observe that if $u$ is a positive radial solution of the form above, then it is necessary that $C_1 > 0$. Hence, from the above asymptotic behavior for the modified Bessel functions, we deduce that

$$u \sim C_1 u_1(r).$$

Now we estimate the constant $C_i$ with $i = 2, 4$. Clearly,

$$K_{\frac{n}{2}-1}(r) \sim_{r \searrow 0} \begin{cases} \log \dfrac{1}{r} & \text{if } n = 2, \\ \dfrac{1}{2}\Gamma\Big(\dfrac{n}{2} - 1\Big)\Big(\dfrac{r}{2}\Big)^{-\frac{n}{2}+1} & \text{if } n > 2, \end{cases}$$

and

$$Y_{\frac{n}{2}-1}(r) \sim_{r \searrow 0} \begin{cases} -\dfrac{2}{\pi} \log \dfrac{1}{r} & \text{if } n = 2, \\ -\dfrac{1}{\pi}\Gamma\Big(\dfrac{n}{2} - 1\Big)\Big(\dfrac{r}{2}\Big)^{-\frac{n}{2}+1} & \text{if } n > 2. \end{cases}$$

see, e.g., [AS64, Sections 9.1 and 9.6]. Therefore, by working around zero and because $u(0)$ is finite, we deduce that

$$C_2 = \frac{2}{\pi} C_4,$$

namely

$$u(r) = C_1 u_1(r) + \frac{2C_4}{\pi} u_2(r) + C_3 u_3(r) + C_4 u_4(r).$$

Taking the Laplacian both sides to get

$$\Delta u(r) = C_1 u_1(r) + \frac{2C_4}{\pi} u_2(r) - C_3 u_3(r) - C_4 u_4(r).$$

Therefore, by working around zero and because $\Delta u(0)$ is finite, we deduce that

$$\frac{2C_4}{\pi} = -\frac{2}{\pi} C_4.$$

Hence $C_4 = 0$, which then implies $C_2 = 0$. Thus, we simply have

$$u(r) = C_1 u_1(r) + C_3 u_3(r).$$

Now taking the Laplacian both sides to get

$$\Delta u(r) = C_1 u_1(r) - C_3 u_3(r).$$

Hence

$$u(r) = C_1 u_1(r) + \frac{u(r) - \Delta u(r)}{2}.$$

In particular, we know that

$$u(0) = C_1 u_1(0) + \frac{u(0) - \Delta u(0)}{2},$$

namely

$$C_1 u_1(0) = \frac{u(0) + \Delta u(0)}{2}$$

From this we are able to compute $C_1$ once we know $u_1(0)$. For $u_1$, we clearly have

$$\begin{aligned} u_1(r) &= r^{-\frac{n}{2}+1} I_{\frac{n}{2}-1}(r) \\ &= r^{-\frac{n}{2}+1} \sum_{l=0}^{+\infty} \frac{1}{l!\,\Gamma(l+n/2)} \Big(\frac{r}{2}\Big)^{2l+n/2-1} \\ &= \frac{2^{-\frac{n}{2}+1}}{\Gamma(n/2)} + 2^{-\frac{n}{2}+1} \sum_{l=1}^{+\infty} \frac{1}{l!\,\Gamma(l+n/2)} \Big(\frac{r}{2}\Big)^{2l}. \end{aligned}$$

Hence we obtain

$$u_1(0) = \frac{2^{-\frac{n}{2}+1}}{\Gamma(n/2)}.$$

Furthermore,

$$u_1(r) = r^{-\frac{n}{2}+1} \frac{e^r}{\sqrt{2\pi r}} \Big[1 - O(r^{-1})_{r \nearrow +\infty}\Big] \sim \frac{1}{\sqrt{2\pi}} r^{-\frac{n-1}{2}} e^r.$$

Hence

$$u(r) \sim \frac{u(0) + \Delta u(0)}{2} \Gamma\Big(\frac{n}{2}\Big) 2^{\frac{n}{2}-1} \frac{1}{\sqrt{2\pi}} r^{-\frac{n-1}{2}} e^r.$$

This is exactly the asymptotic behavior of $u$ as wanted. □

As can be seen from the above proof, the two regular solutions $u_1$ and $u_3$ are written in terms of the Bessel and modified Bessel functions.

*Remark* 3.1. In other words, the general solution to $\Delta^2 u(r) = u(r)$ in $\mathbf{R}^n$, which depends only on two parameters, is actually a linear combination of $I$ and $J$.

3.3. **The common cases**: $\alpha \in (-1, 1)$ **with** $n \geq 2$, **and** $\alpha \in (-1/3, 1)$ **if** $n = 1$. Now we consider the case $-1 < \alpha < 1$. The remaining cases $\alpha \leq -1$ will be treated later.

**Proposition 3.2**. *Assume that $\alpha \in (-1,1)$ if $n \geq 2$, or $\alpha \in (-1/3,1)$ if $n = 1$. Let $u$ be a positive radial solution to the problem* (1.1) *in $\mathbf{R}^n$, then $u$ has the following asymptotic*

*behavior*

$$u \sim C(n,\alpha)^{-\frac{1}{1-\alpha}} r^{\frac{4}{1-\alpha}}.$$

*Proof.* The proof consists of two cases:

**Case 1**. Suppose that $\alpha \in [0,1)$. Recall from Lemma 3.2 the existence of $C_0 > 0$ and $r_0 > 0$ such that $u(r) \geq C_0 r^2$ and $\Delta u > 0$ for all $r \geq r_0$. In particular, $u$ and $\Delta u$ are monotone increasing in $(r_0, +\infty)$. Hence, for $r > r_0$ we easily get

$$\Delta u(2r) = \Delta u(r) + \int_r^{2r} \frac{1}{|\mathbb{S}^{n-1}|t^{n-1}} \int_{B_t} u^\alpha(x) dx dt \geq c_1 r^2 u^\alpha(r)$$

for some constant $c_1 > 0$. We also get

$$u(4r) = u(2r) + \int_{2r}^{4r} \frac{1}{|\mathbb{S}^{n-1}|t^{n-1}} \int_{B_t} \Delta u(x) dx dt \geq c_2 r^2 \Delta u(2r) \geq c_1 c_2 r^4 u^\alpha(r)$$

for some constant $c_2 > 0$. From this, if we denote

$$U(r) = r^{-\frac{4}{1-\alpha}} u(r) \quad \text{with } r > 0,$$

then the above estimate can be rewritten as follows

$$U(4r) \geq \underbrace{4^{-\frac{4}{1-\alpha}} c_1 c_2}_{=c} U^\alpha(r) \quad \text{with } r > r_0. \tag{3.2}$$

Now we show that

$$\liminf_{r \nearrow \infty} U(r) > 0.$$

Obviously, $U(r)$ is bounded from below in $B_{4r_0} \setminus B_{r_0}$. By using (3.2) we easily obtain

$$U(4^k r) \geq c^{1+\alpha+\cdots+\alpha^{k-1}} U^{\alpha^k}(r)$$

for any $r > r_0$ and any integer $k \geq 1$. As $\alpha \in [0,1)$, it is crucial to have

$$\lim_{k \nearrow \infty} c^{1+\alpha+\cdots+\alpha^{k-1}} = c^{\frac{1}{1-\alpha}}$$

and

$$\lim_{k \nearrow \infty} \alpha^k = 0.$$

This and the boundedness of $U$ in the ring $B_{4r_0} \setminus B_{r_0}$ allow us to conclude that $U$ cannot be too small near infinity and this is enough to conclude that $\liminf_{r \nearrow \infty} U(r) > 0$. An immediate consequence is the following

$$u(r) \geq C_2 r^{\frac{4}{1-\alpha}} \tag{3.3}$$

for some $C_2 > 0$ and for any $r \geq 0$. (In the earlier version [**NNP18**], the preceding estimate was proved by using a pointwise comparison principle for higher-order elliptic equations, see [**FF16**, Proposition A.2]. Here we offer an alternative argument.) Since $\alpha \geq 0$, such a lower bound for $u$ implies that

$$\Delta^2 u(r) \geq C_3 r^{\frac{4\alpha}{1-\alpha}}$$

for any $r \geq 0$. Integrating the above differential inequality twice to get

$$\Delta u(r) \geq C_4 r^{2+\frac{4\alpha}{1-\alpha}},$$

for some $C_4 > 0$ and for any $r \geq r_1 \gg 1$. Consequently, $u$ is increasing in $(r_1, +\infty)$. Then for any $r > r_1$ we have from Lemma 3.1 the following

$$\Delta u(r) = \Delta u(r_1) + \int_0^{r_1} t^{n-1} u^\alpha(t) dt \int_{r_1}^r s^{1-n} ds + \int_{r_1}^r s^{1-n} \Big( \int_{r_1}^s t^{n-1} u^\alpha(t) dt \Big) ds$$

$$\leq \Delta u(r_1) + \int_0^{r_1} t^{n-1} u^\alpha(t)dt \int_{r_1}^r s^{1-n} ds + u^\alpha(r) \int_{r_1}^r s^{1-n}\Big(\int_{r_1}^s t^{n-1} dt\Big) ds.$$

Since

$$\int_{r_1}^r s^{1-n}\Big(\int_{r_1}^s t^{n-1} dt\Big) ds = \begin{cases} \dfrac{(r-r_1)^2}{2} & \text{if } n=1, \\ \dfrac{(r-r_1)^2}{2n} - \dfrac{r_1^n\big(r_1^{2-n} - r^{2-n}\big)}{n(n-2)} & \text{if } n \geq 2, \end{cases}$$

and

$$\int_{r_1}^r s^{1-n} ds = \begin{cases} r - r_1 & \text{if } n=1, \\ \log\Big(\dfrac{r}{r_1}\Big) & \text{if } n=2, \\ (n-2)^{-1}\big(r_1^{2-n} - r^{2-n}\big) & \text{if } n \geq 3, \end{cases}$$

we deduce that

$$\Delta u(r) \leq C_5 + C_5 r + C_5 r^2 u^\alpha(r)$$

for some $C_5 > 0$ and for any $r \geq r_1$. From this and (3.3), it follows that

$$\Delta u(r) \leq C_6 r^2 u^\alpha(r)$$

for some $C_6 > 0$ and for any $r \geq r_2 \gg r_1$. Keep in mind that $\Delta u$ is increasing in $(0, +\infty)$. Therefore, we repeat the above argument to get

$$u(r) \leq \Delta u(r_1) + \int_0^{r_1} t^{n-1} \Delta u(t)dt \int_{r_1}^r s^{1-n} ds + \Delta u(r) \int_{r_1}^r s^{1-n}\Big(\int_{r_1}^s t^{n-1} dt\Big) ds.$$

From this, it is not hard to check that

$$u(r) \leq C_7 r^2 \Delta u(r)$$

for some $C_7 > 0$ and for any $r \geq r_3 \gg r_2$. Hence, we have just shown that

$$u(r) \leq C_6 C_7 r^4 u^\alpha(r)$$

for any $r \geq r_3$. This implies that $u$ has the following upper bound

$$u(r) \leq C_8 r^{\frac{4}{1-\alpha}} \tag{3.4}$$

for some $C_8 > 0$ and for any $r \geq r_3$. From (3.3) and (3.4), we have

$$0 < \liminf_{r \nearrow +\infty} r^{-\frac{4}{1-\alpha}} u(r) =: a_{\inf} \leq a_{\sup} := \limsup_{r \nearrow +\infty} r^{-\frac{4}{1-\alpha}} u(r) < +\infty.$$

For any $\epsilon \in (0, a_{\inf})$, there exists $R(\epsilon) > 0$ such that

$$a_{\inf} - \epsilon < r^{-\frac{4}{1-\alpha}} u(r) < a_{\sup} + \epsilon$$

for any $r \geq R(\epsilon)$. Under the condition $\alpha > 0$ and for any $r \geq R(\epsilon)$, we deduce from the second identity in Lemma 3.1 the following

$$\begin{aligned} \Delta u(r) \leq& \Delta u(R(\epsilon)) + \int_0^{R(\epsilon)} t^{n-1} u^\alpha(t)dt \int_{R(\epsilon)}^r s^{1-n} ds \\ &+ (a_{\sup} + \epsilon)^\alpha \int_{R(\epsilon)}^r s^{1-n}\Big(\int_{R(\epsilon)}^s t^{n-1+\frac{4\alpha}{1-\alpha}} dt\Big) ds \\ \leq& O(1) r + \frac{(a_{\sup}+\epsilon)^\alpha}{(n + \frac{4\alpha}{1-\alpha})(2 + \frac{4\alpha}{1-\alpha})} r^{2+\frac{4\alpha}{1-\alpha}} \end{aligned}$$

with error $O(1)$ depending only on $\epsilon$. From this, by integrating twice, we arrive at

$$u(r) \leq O(1) r^3 + \frac{(a_{\sup}+\epsilon)^\alpha}{\prod_{l=1}^2 (n + \frac{4\alpha}{1-\alpha} + 2l - 2)(2l + \frac{4\alpha}{1-\alpha})} r^{\frac{4}{1-\alpha}}$$

for any $r \geq R(\epsilon)$. Observe that $4/(1-\alpha) > 3$. Hence, simply dividing both side by $r^{4/(1-\alpha)}$ and letting $r \nearrow +\infty$, we get

$$a_{\sup} \leq \frac{(a_{\sup}+\epsilon)^\alpha}{\prod_{l=1}^2 (n+\frac{4\alpha}{1-\alpha}+2l-2)(2l+\frac{4\alpha}{1-\alpha})}$$

for any $\epsilon \in (0, a_{\inf})$. Now letting $\epsilon \to 0$ to get

$$a_{\sup} \leq \Big( \prod_{l=1}^2 \Big(n+\frac{4\alpha}{1-\alpha}+2l-2\Big)\Big(2l+\frac{4\alpha}{1-\alpha}\Big)\Big)^{-\frac{1}{1-\alpha}}. \tag{3.5}$$

By the same argument, we get

$$\begin{aligned}
\Delta u(r) \geq & \Delta u(R(\epsilon)) + \int_0^{R(\epsilon)} t^{n-1} u^\alpha(t) dt \int_{R(\epsilon)}^r s^{1-n} ds \\
& + (a_{\inf}-\epsilon)^\alpha \int_{R(\epsilon)}^r s^{1-n} \Big( \int_{R(\epsilon)}^s t^{n-1+\frac{4\alpha}{1-\alpha}} dt \Big) ds \\
= & O(1) r + \frac{(a_{\inf}-\epsilon)^\alpha}{(n+\frac{4\alpha}{1-\alpha})(2+\frac{4\alpha}{1-\alpha})} \Big[ r^{2+\frac{4\alpha}{1-\alpha}} - R(\epsilon)^{2+\frac{4\alpha}{1-\alpha}} \Big] \\
& - \frac{(a_{\inf}-\epsilon)^\alpha}{(n+\frac{4\alpha}{1-\alpha})(2-n)} R(\epsilon)^{n+\frac{4\alpha}{1-\alpha}} \Big[ r^{2-n} - R(\epsilon)^{2-n} \Big] \\
\geq & O(1) r + \frac{(a_{\inf}-\epsilon)^\alpha}{(n+\frac{4\alpha}{1-\alpha})(2+\frac{4\alpha}{1-\alpha})} r^{2+\frac{4\alpha}{1-\alpha}} - C(n,\alpha,\epsilon)
\end{aligned}$$

with error $O(1)$ depending only on $\epsilon$. Here the constant $C > 0$ depends only on $n$, $\alpha$, and $\epsilon$. From this, integrating twice gives

$$u(r) \geq O(1) r^3 + \frac{(a_{\inf}-\epsilon)^\alpha}{\prod_{l=1}^2 (n+\frac{4\alpha}{1-\alpha}+2l-2)(2l+\frac{4\alpha}{1-\alpha})} r^{\frac{4}{1-\alpha}}$$

for any $r \geq R(\epsilon)$. Hence as above we can easily prove that

$$a_{\inf} \geq \Big( \prod_{l=1}^2 \Big(n+\frac{4\alpha}{1-\alpha}+2l-2\Big)\Big(2l+\frac{4\alpha}{1-\alpha}\Big)\Big)^{-\frac{1}{1-\alpha}}. \tag{3.6}$$

Combining (3.5) and (3.6), we get

$$\lim_{r \nearrow +\infty} r^{-\frac{4}{1-\alpha}} u(r) = \Big( \prod_{l=1}^2 \Big(n+\frac{4\alpha}{1-\alpha}+2l-2\Big)\Big(2l+\frac{4\alpha}{1-\alpha}\Big)\Big)^{-\frac{1}{1-\alpha}}$$

as indicated. This completes the first case.

**Case 2**. Suppose that $\alpha < 0$. Then in this scenario, we shall consider either $\alpha \in (-1,0)$ if $n \geq 2$ or $\alpha \in (-1/3, 0)$ if $n = 1$. Note by Lemma 3.2 that $\lim_{r\nearrow+\infty} \Delta u(r) > 0$. This implies that

$$\int_0^R s^{n-1} \Delta u(s) ds > 0$$

for $R$ large. Hence, there exists $R_0 > 0$ such that

$$u'(R) = R^{1-n} \int_0^R s^{n-1} \Delta u(s) ds > 0$$

for any $R \geq R_0$. In other words, $u$ is increasing on $[R_0, +\infty)$. We now make use of the rescaled test-function argument. Indeed, let $\psi = \psi(r)$ be a smooth radial cut-off function satisfying $0 \le \psi \le 1$ and

$$\psi(r) = \begin{cases} 0 & \text{if } r \in [0,1] \cup [4,+\infty], \\ 1 & \text{if } 2 \le r \le 3. \end{cases}$$

For any $R \ge R_0$, let $\phi_R(r) = \psi(r/R)$. Then we have

$$\begin{aligned} \int_0^{+\infty} u^\alpha \phi_R s^{n-1} ds &= \int_0^{+\infty} \Delta^2 u \, \phi_R s^{n-1} ds \\ &= \int_0^{+\infty} u \, \Delta^2 \phi_R s^{n-1} ds \le C_1^{-1} R^{-4} \int_R^{4R} u \, s^{n-1} ds \end{aligned}$$

for some $C_1 > 0$. Hence,

$$R^{-4} \int_R^{4R} u \, s^{n-1} ds \ge C_1 \int_0^{+\infty} u^\alpha \phi_R s^{n-1} ds \ge C_1 \int_{2R}^{3R} u^\alpha s^{n-1} ds.$$

Using $\alpha < 0$ and the monotonicity of $u$ on $[R_0, +\infty)$, we deduce that

$$R^{-4} R^n u(4R) \ge C_1 u^\alpha(4R) R^n$$

for any $R \ge R_0$. Consequently, we obtain the lower bound

$$u(R) \ge C_1 R^{\frac{4}{1-\alpha}} \tag{3.7}$$

for any $R \ge R_1 = 4R_0$. By the assumption that either $\alpha \in (-1,0)$ if $n \ge 2$ or $\alpha \in (-1/3, 0)$ if $n = 1$, we can always have that $n + 4\alpha/(1-\alpha) > 0$. Thanks to the monotonicity of $u$ and $\Delta u$ for $R$ large and the fact that $\lim_{r \nearrow +\infty} \Delta u(R) > 0$, we follow the same argument used in the previous case to obtain

$$u(R) \le C_2 R^2 \Delta u(R) \tag{3.8}$$

for some $C_2 > 0$ and for any $R \ge R_2 \gg R_1$. On the other hand, for any $R \ge R_2$ we can estimate

$$\begin{aligned} \Delta u(R) = \Delta u(R_2) &+ \int_0^{R_2} t^{n-1} u^\alpha(t) dt \int_{R_2}^R s^{1-n} ds \\ &+ \int_{R_2}^R s^{1-n} \Big( \int_{R_2}^s t^{n-1} u^\alpha(t) dt \Big) ds \\ \le \Delta u(R_2) &+ \int_0^{R_2} t^{n-1} u^\alpha(t) dt \int_{R_2}^R s^{1-n} ds \\ &+ C_2^\alpha \int_{R_2}^R s^{1-n} \Big( \int_{R_2}^s t^{n-1+\frac{4\alpha}{1-\alpha}} dt \Big) ds, \end{aligned}$$

thanks to (3.7) and $\alpha < 0$. By the previous inequality and the fact $n + 4\alpha/(1-\alpha) > 0$, simply considering either $n = 1$ or $n \ge 2$ separately, we can find a constant $C_3 > 0$ and $R_3 \gg R_2$ such that

$$\Delta u(R) \le C_3 R^{2+\frac{4\alpha}{1-\alpha}} \tag{3.9}$$

for any $R \ge R_3$. Combining (3.8) and (3.9), we obtain the upper bound

$$u(R) \le C_4 R^{\frac{4}{1-\alpha}} \tag{3.10}$$

for some $C_4 > 0$ and $R \ge R_3$.

Once we can bound $u$ from above and below as shown in (3.7) and (3.10), we can repeat the argument used in Case 1 to obtain the desired limit. Indeed, we let $a_{\inf}$ and $a_{\sup}$ be the following

$$0 < \liminf_{R \to +\infty} R^{-\frac{4}{1-\alpha}} u(R) =: a_{\inf} \le a_{\sup} := \limsup_{R \to +\infty} R^{-\frac{4}{1-\alpha}} u(R) < +\infty.$$

For any $\epsilon \in (0, a_{\inf})$, there exists $R(\epsilon) > 0$ such that

$$a_{\inf} - \epsilon < r^{-\frac{4}{1-\alpha}} u(R) < a_{\sup} + \epsilon$$

for any $R \geq R(\epsilon)$. Under the condition $\alpha < 0$, for any $R \geq R(\epsilon)$, we have

$$\begin{aligned}\Delta u(R) \leq & \Delta(R(\epsilon)) + \int_0^{R(\epsilon)} t^{n-1} u^{\alpha(t)} dt \int_{R(\epsilon)}^R s^{1-n} ds \\ & + (a_{\inf} - \epsilon)^\alpha \int_{R(\epsilon)}^R s^{1-n} \Big( \int_{R(\epsilon)}^s t^{n-1+\frac{4\alpha}{1-\alpha}} dt \Big) ds \\ = & \begin{cases} O(1)R + \big(n + \dfrac{4\alpha}{1-\alpha}\big)^{-1} \big(2 + \dfrac{4\alpha}{1-\alpha}\big)^{-1} (a_{\inf} - \epsilon)^\alpha R^{2+\frac{4\alpha}{1-\alpha}} & \text{if } n = 1, \\ O(1) \log R + \big(n + \dfrac{4\alpha}{1-\alpha}\big)^{-1} \big(2 + \dfrac{4\alpha}{1-\alpha}\big)^{-1} (a_{\inf} - \epsilon)^\alpha R^{2+\frac{4\alpha}{1-\alpha}} & \text{if } n \geq 2. \end{cases}\end{aligned}$$

From this, by integrating by parts, we arrive at

$$u(R) \leq \begin{cases} O(1)R^3 + \big(n + \dfrac{4\alpha}{1-\alpha}\big)^{-1} \big(2 + \dfrac{4\alpha}{1-\alpha}\big)^{-1} (a_{\inf} - \epsilon)^\alpha R^{\frac{4}{1-\alpha}} & \text{if } n = 1, \\ O(1)R^2 \log R + \big(n + \dfrac{4\alpha}{1-\alpha}\big)^{-1} \big(2 + \dfrac{4\alpha}{1-\alpha}\big)^{-1} (a_{\inf} - \epsilon)^\alpha R^{\frac{4}{1-\alpha}} & \text{if } n \geq 2, \end{cases}$$

for any $R \geq R(\epsilon)$. Keep in mind that $4/(1-\alpha) > 3$ if $n = 1$ and $4/(1-\alpha) > 2$ if $n \geq 2$. Hence, dividing both side by $R^{4/(1-\alpha)}$ and letting $R \to +\infty$, we get

$$a_{\sup} \leq \frac{(a_{\inf} - \epsilon)^\alpha}{\prod_{l=1}^2 (n + \frac{4\alpha}{1-\alpha} + 2l - 2)(2l + \frac{4\alpha}{1-\alpha})}$$

for any $\epsilon \in (0, a_{\inf})$. Letting $\epsilon \to 0$, we get

$$a_{\sup} a_{\inf}^{-\alpha} \leq \frac{1}{\prod_{l=1}^2 (n + \frac{4\alpha}{1-\alpha} + 2l - 2)(2l + \frac{4\alpha}{1-\alpha})}. \tag{3.11}$$

By the same argument, we can easily prove that

$$a_{\inf} a_{\sup}^{-\alpha} \geq \frac{1}{\prod_{l=1}^2 (n + \frac{4\alpha}{1-\alpha} + 2l - 2)(2l + \frac{4\alpha}{1-\alpha})}. \tag{3.12}$$

Combining (3.11) and (3.12), we get

$$\Big(\frac{a_{\sup}}{a_{\inf}}\Big)^{1+\alpha} \leq 1.$$

Since $1 + \alpha > 0$, we clearly have $a_{\sup} \leq a_{\inf}$. From this, we must have $a_{\sup} = a_{\inf}$ and therefore

$$\lim_{r \nearrow +\infty} r^{-\frac{4}{1-\alpha}} u(r) = \Big( \prod_{l=1}^2 \Big(n + \frac{4\alpha}{1-\alpha} + 2l - 2\Big)\Big(2l + \frac{4\alpha}{1-\alpha}\Big)\Big)^{-\frac{1}{1-\alpha}} = C(n, \alpha)^{-\frac{1}{1-\alpha}}$$

as claimed. □

Using a well-know comparison principle for polyharmonic operator, we can also prove (3.3) quickly. We also note that from the above proof by using the upper and the lower bounds for $\Delta u$ we can also obtain the asymptotic behavior of $\Delta u$ near infinity.

3.4. **The case $\alpha \leq -1$ with $n \geq 3$.** We now consider the case $\alpha \leq -1$ covered in Theorem 1. This case is split into two sub-cases corresponding to either $\alpha = -1$ or not. First we consider the case $\alpha = -1$.

**Proposition 3.3.** *Assume that $\alpha = -1$ and $n \geq 3$. Let $u$ be a positive radial solution to the problem (1.1) in $\mathbf{R}^n$, then $u$ has the following asymptotic behavior*

$$u(r) \sim \frac{1}{\sqrt{n(n-2)}} r^2 \sqrt{\log r}.$$

*Proof.* It follows from Lemma 3.2 that

$$\gamma = \lim_{r \nearrow +\infty} \Delta u(r) > 0,$$

which could be infinity. Suppose that $\gamma < +\infty$. By Lemma 3.2 we know that $u \sim (\gamma/(2n))r^2$ at infinity. Using this behavior and Lemma 3.1, we can easily bound $\Delta u$ from below as shown below

$$\begin{aligned}
\Delta u(r) &= \Delta u(R) + \int_0^R t^{n-1} u^{-1} dt \int_R^r s^{1-n} ds + \int_R^r s^{1-n} \Big( \int_R^s t^{n-1} u^{-1}(t) dt \Big) ds \\
&\geq \Delta u(R) + C^{-1} \int_R^r s^{1-n} \int_R^s t^{n-3} dt ds \\
&= \Delta u(R) + \frac{R^{n-2}}{(n-2)^2 C} \frac{1}{r^{n-2}} - \frac{(n-2)\log R + 1}{(n-2)^2 C} + \frac{\log r}{(n-2)C}.
\end{aligned}$$

Clearly, by sending $r \nearrow +\infty$ we obtain a contradiction, thanks to $\gamma < +\infty$. Thus, we must have

$$\lim_{r \nearrow +\infty} \Delta u(r) = +\infty.$$

Using this, we can evaluate $u^{-1}(r)$ as follows:

$$\lim_{r \nearrow +\infty} \frac{r^2}{u(r)} = 2 \lim_{r \nearrow +\infty} \frac{r^n}{r^{n-1} u'(r)} = 2n \lim_{r \nearrow +\infty} \frac{r^{n-1}}{(r^{n-1} u'(r))'} = 0.$$

Thus, this gives $u^{-1}(r) = o(r^{-2})$. Next, using integration by parts, we get from $\Delta^2 u = u^{-1}$ the following

$$\Delta u(r) = \Delta u(0) - \frac{r^{2-n}}{n-2} \int_0^r s^{n-1} u^{-1} ds + \frac{1}{n-2} \int_0^r s u^{-1} ds. \tag{3.13}$$

Now it follows from the l'Hôpital rule and $u^{-1}(r) = o(r^{-2})$ that

$$\lim_{r \nearrow +\infty} \frac{\int_0^r s^{n-1} u^{-1}(s) ds}{r^{n-2}} = \lim_{r \nearrow +\infty} \frac{r^{n-1} u^{-1}(r)}{(n-2) r^{n-3}} = \frac{1}{n-2} \lim_{r \nearrow +\infty} \frac{r^2}{u(r)} = 0.$$

Hence, (3.13) gives

$$\lim_{r \nearrow +\infty} \frac{\Delta u(r)}{F(r)} = \frac{1}{n-2}, \tag{3.14}$$

where

$$F(r) = \int_0^r s u^{-1} ds.$$

Using the l'Hôpital rule and noting that $u^{-1}(r) = o(r^{-2})$, we have

$$\begin{aligned}
\lim_{r \nearrow +\infty} \frac{u'(r)}{rF(r)} &= \lim_{r \nearrow +\infty} \frac{r^{n-1} u'(r)}{r^n F(r)} \\
&= \lim_{r \nearrow +\infty} \frac{(r^{n-1} u'(r))'}{n r^{n-1} F(r) + r^{n+1} u^{-1}(r)} \\
&= \lim_{r \nearrow +\infty} \frac{\Delta u(r)}{nF(r) + r^2 u^{-1}(r)} \\
&= \lim_{r \nearrow +\infty} \frac{\Delta u(r)}{nF(r)}
\end{aligned}$$

$$= \frac{1}{n(n-2)},$$

which helps us to conclude that

$$\lim_{r\nearrow+\infty}\frac{u(r)}{r^2F(r)} = \lim_{r\nearrow+\infty}\frac{u'(r)}{2rF(r)+r^3u^{-1}(r)} = \lim_{r\nearrow+\infty}\frac{u'(r)}{2rF(r)} = \frac{1}{2n(n-2)}. \tag{3.15}$$

Replacing $u(r) = r(F'(r))^{-1}$, we deduce from (3.15) that

$$\lim_{r\nearrow+\infty}\frac{1}{rF'(r)F(r)} = \frac{1}{2n(n-2)}.$$

It follows from the l'Hôpital rule that

$$\lim_{r\nearrow+\infty}\frac{2\log r}{F^2(r)} = \frac{1}{2n(n-2)}. \tag{3.16}$$

The proof follows by combining (3.15) and (3.16). □

From the above proof, we can also obtain a precise asymptotic behavior of $\Delta u$ near infinity. Indeed, if we combine (3.14) and (3.16), the we easily get

$$\Delta u(r) \sim 2\sqrt{\frac{n}{n-2}}(\log r)^{1/2}.$$

We now consider the case $\alpha < -1$.

**Proposition 3.4.** *Assume that $n \geq 3$ and $\alpha < -1$. Let $u$ be a positive radial solution to the problem (1.1) in $\mathbf{R}^n$, then $u$ has the following asymptotic behavior*

$$u \sim \frac{1}{2n}\Big(\Delta u(0) + \frac{1}{n-2}\int_0^{+\infty} tu^\alpha(t)dt\Big)r^2.$$

*Proof.* It follows from Lemma 3.2 that $\lim_{r\nearrow+\infty}\Delta u(r) > 0$ and that $u(r) \geq Cr^2$ for any $r \geq R_0$. Hence, thanks to $\alpha < -1$, the integral

$$\int_0^{+\infty} s^{1-n}\int_0^s t^{n-1}u^\alpha(t)dtds$$

exists. From this, using the representation formula

$$\Delta u(r) = \Delta u(0) + \int_0^r s^{1-n}\Big(\int_0^s t^{n-1}u^\alpha(t)dt\Big)ds$$

we deduce that $\Delta u(r)$ has the finite limit as $r \nearrow +\infty$ and, in addition, there holds

$$\begin{aligned}\lim_{r\nearrow+\infty}\Delta u(r) &= \Delta u(0) + \int_0^{+\infty} s^{1-n}\Big(\int_0^s t^{n-1}u^\alpha(t)dt\Big)ds\\ &= \Delta u(0) + \frac{1}{n-2}\int_0^{+\infty} tu^\alpha(t)dt.\end{aligned}$$

The proof is complete by making use of Lemma 3.3 . □

3.5. **The case $\alpha \leq -1$ with $n = 2$.** This subsection is devoted to proofs of the remaining cases in Theorem 2 indicated in Table 1 since the cases $\alpha = 1$ and $\alpha \in (-1, 1)$ were already proved in Propositions 3.1 and 3.2 respectively. As always, we consider the cases $\alpha = -1$ and $\alpha < -1$ separately.

**Proposition 3.5.** *Assume that $\alpha = -1$ and $n = 2$. Let $u$ be a positive radial solution to the problem (1.1) in $\mathbf{R}^n$, then $u$ has the following asymptotic behavior*

$$u \sim \frac{1}{\sqrt{2}} r^2 \log r \sqrt{\log\log r}.$$

*Proof.* In view of Lemma 3.2, we deduce that

$$\gamma = \lim_{r \nearrow +\infty} \Delta u(r) > 0.$$

As in the proof of Proposition 3.3, if $\gamma < +\infty$, then there exist two numbers $C > 0$ and $R > 10$ such that

$$u(r) \le C r^2$$

for any $r \ge R$. Let $r \nearrow +\infty$ in the inequality

$$\begin{aligned} \Delta u(r) &= \Delta u(0) + \int_0^r s^{-1}\Big(\int_0^s t u^{-1}(t) dt\Big) ds \\ &\ge \Delta u(0) + \int_1^r s^{-1}\Big(\int_0^1 t u^{-1}(t) dt\Big) ds, \end{aligned}$$

we know that $\Delta u(r) \to +\infty$ as $r \nearrow +\infty$, which contradicts $\gamma < +\infty$. Hence, this proves $\gamma = +\infty$. Consequently, we have $u(r) \ge C_1 r^2$ for some $C_1 > 0$ and for any $r \ge R$. For simplicity, we set

$$G(r) = \int_0^r t u^{-1}(t) dt.$$

Keep in mind that $G(r)\log r \nearrow +\infty$ as $r \nearrow +\infty$. Thus, using integration by parts, we get

$$\begin{aligned} \Delta u(r) &= \Delta u(0) + \int_0^r s^{-1} G(s) ds \\ &= \Delta u(0) + G(r)\log r - \int_0^r t u^{-1}(t) \log t dt \\ &\le \Delta u(0) - \int_0^1 t \log t u^{-1}(t) dt + G(r)\log r \\ &\le C_2 G(r) \log r \end{aligned}$$

for some $C_2 > 0$ and for any $r \ge R$. Using this and the monotonicity of $\Delta u$, we can bound $u$ from above as follows

$$\begin{aligned} u(r) &= u(R) + G(R)\int_R^r s^{-1} ds + \int_R^r s^{-1}\Big(\int_R^s t \Delta u(t) dt\Big) ds \\ &\le u(R) + G(R)\int_R^r s^{-1} ds + \Delta u(r) \int_R^r s^{-1}\Big(\int_R^s t dt\Big) ds \\ &\le C_3 G(r) r^2 \log r \end{aligned}$$

for some $C_3 > 0$ and for $r \ge R_1 \gg R$. Consequently,

$$r u(r)^{-1} G(r) \ge C_3 r^{-1} (\log r)^{-1},$$

for $r \ge R_1$. Integrating this inequality over $(R_1, r)$, we obtain

$$G(r)^2 - G(R_1)^2 \ge C_3\big(\log\log r - \log\log R_1\big),$$

which implies that

$$G(r) \ge C_4 \sqrt{\log\log r},$$

for some $C_4 > 0$ and for all $r \geq R_2$ for some $R_2 \gg R_1$. Note that

$$\Delta u(r) = \Delta u(0) + \int_0^r s^{-1} G(s) ds.$$

Hence for $r \geq R_2^2$, we can estimate

$$\Delta u(r) \geq \int_{\sqrt{r}}^r s^{-1} G(s) ds \geq C_4 \sqrt{\log\log\sqrt{r}} \log\sqrt{r},$$

which gives

$$\Delta u(r) \geq C_5 \log r \sqrt{\log\log r}$$

for some $C_5 > 0$ and for any $r \geq R_3 \gg R_2^2$. Making use of this inequality in the integral expression of $u$, that is

$$u(r) = u(0) + \int_0^r s^{-1} \int_0^s t \Delta u(t) dt ds,$$

and performing a similar argument as above, we eventually get

$$u(r) \geq C_6 r^2 \log r \sqrt{\log\log r}$$

for some $C_6 > 0$ and for any $r \geq R_4 \gg R_3$. This and the equation $\Delta^2 u = u^{-1}$ imply that

$$\Delta^2 u(r) \leq C_6^{-1} r^{-2} (\log r)^{-1} (\log\log r)^{-1/2}$$

for $r \geq R_4$. Since we are in $\mathbf{R}^2$, integrating the above differential inequality to get

$$(\Delta u)'(r) - (\Delta u)'(R_4) \leq 2C_6^{-1} (\sqrt{\log\log r} - \sqrt{\log\log R_4}),$$

which implies that

$$(\Delta u)'(r) \leq C_7 \sqrt{\log\log r}$$

for some $C_7 > 0$ and for any $r \geq R_5 \gg R_4$. Continuing this process, we arrive at

$$\Delta u(r) \leq C_8 \log r \sqrt{\log\log r}$$

for some $C_8 > 0$ and for any $r \geq R_6 \gg R_5$. Simply repeating the above argument, we eventually get

$$u(r) \leq C_9 r^2 \log r \sqrt{\log\log r}$$

for some $C_9 > 0$ and for any $r \geq R_7 \gg R_6$. Thus, we have already shown that

$$C_6 r^2 \log r \sqrt{\log\log r} \leq u(r) \leq C_9 r^2 \log r \sqrt{\log\log r} \tag{3.17}$$

for $r$ large. Applying the first identity from Lemma 3.1 and integrating by parts over $[R_7, r]$ to get

$$\Delta u(r) = \Delta u(R_7) + \Big( G(r) \log r - G(R_7) \log R_7 \Big) - \int_{R_7}^r t u^{-1}(t) \log t dt. \tag{3.18}$$

Using the lower bound for $G$, we deduce that

$$G(r) \log r \geq C_4 \sqrt{\log\log r} \log r.$$

We now use (3.17) to bound

$$t u^{-1}(t) \log t \leq C_6^{-1} t^{-1} (\sqrt{\log\log t})^{-1}$$

for all $t \geq R_7$. From this we obtain

$$\int_{R_7}^r t u^{-1}(t) \log t \, dt \leq \int_{r/2}^r t u^{-1}(t) \log t \, dt \leq C_6^{-1} \log\log(r/2))^{-1/2} \log r$$

for any $r \geq 2R_7$. Hence, dividing both sides of (3.18) by $G(r)\log r$ and sending $r$ to infinity to get

$$\lim_{r \nearrow +\infty} \frac{\Delta u(r)}{G(r)\log r} = 1.$$

Using the l'Hôpital rule and noting that $r^2 = o(u(r)G(r))$, we have

$$\lim_{r \nearrow +\infty} \frac{u'(r)}{G(r) r \log r} = \lim_{r \nearrow +\infty} \frac{(ru'(r))'}{r^3 u^{-1}(r) \log r + 2G(r) r \log r + G(r) r} = \frac{1}{2}.$$

Applying the l'Hôpital rule one more time, we deduce from the preceding limit that

$$\lim_{r \nearrow +\infty} \frac{u(r)}{G(r) r^2 \log r} = \frac{1}{4}. \tag{3.19}$$

Next, replacing $u(r) = r(G'(r))^{-1}$, we deduce from (3.19) that

$$\lim_{r \nearrow +\infty} \frac{1}{G'(r)G(r) r \log r} = \frac{1}{4}.$$

Hence,

$$\lim_{r \nearrow +\infty} \frac{2\log(\log r)}{G^2(r)} = \frac{1}{4}. \tag{3.20}$$

The proof follows by putting (3.19) and (3.20) together. □

Next we consider the case $\alpha < -1$. Our result for this case is the following.

**Proposition 3.6.** *Assume that $n = 2$ and $\alpha < -1$. Let $u$ be a positive radial solution to the problem (1.1) in $\mathbf{R}^n$, then $u$ has the following asymptotic behavior*

$$u \sim \Big(\frac{1}{4} \int_0^{+\infty} t u^\alpha(t) dt\Big) r^2 \log r.$$

*Proof.* We recall from Lemma 3.2 that $u(r) \geq C_0 r^2$ for any $r \geq r_0$. This and the fact $\alpha < -1$ allows us to estimate

$$0 < D := \int_0^{+\infty} t u^\alpha(t) dt < +\infty.$$

Now it follows from

$$r(\Delta u)'(r) = \int_0^r t u^\alpha(t) dt$$

that

$$\lim_{r \nearrow +\infty} r(\Delta u)'(r) = D.$$

Hence, applying the L'Hôpital rule gives

$$\lim_{r \nearrow +\infty} \frac{u(r)}{r^2 \log r} = \lim_{r \nearrow +\infty} \frac{\Delta u(r)}{4 \log r} = \lim_{r \nearrow +\infty} \frac{r(\Delta u)'(r)}{4} = \frac{D}{4}$$

as claimed. □

3.6. **The case $\alpha \leq -1/3$ with $n = 1$.** This subsection is devoted to proofs of Theorems Theorem 3 indicated in Table 1. Since the cases $\alpha = 1$ and $\alpha \in (-1/3, 1)$ are proved in Propositions 3.1 and 3.2 respectively. We only give the proof for $\alpha \leq -1/3$ in this sub-section. First we consider the case $\alpha = -1/3$.

**Proposition 3.7.** *Assume that* $n = 1$ *and* $\alpha = -1/3$. *Let* $u$ *be a positive radial solution to the problem* (1.1) *in* $\mathbf{R}^n$, *then* $u$ *has the following asymptotic behavior*

$$u(r) \sim \Big(\frac{2}{9}\Big)^{3/4} r^3 (\log r)^{3/4}.$$

*Proof.* We follow the argument as in the proof of (3.7) to get the following estimate

$$u(R) \geq C_1 R^{\frac{4}{1-\alpha}}$$

for some $C_1 > 0$ and $R$ large enough. Recall that $u$ solves $u^{(4)} = u^{-1/3}$ in $\mathbf{R}$ and $u^{(3)}(0) = 0$. From this, by integration by parts, we get

$$u''(r) = u''(0) + r\int_0^r u^{-1/3}(s)ds - \int_0^r su^{-1/3}(s)ds \leq C_2 r\int_0^r u^{-1/3}(s)ds$$

for some $C_2 > 0$ and for any $r \geq R$. Hence, thanks to the first identity in Lemma 3.1, we obtain

$$\begin{aligned} u(r) &= u(R) + \int_R^r \Big(\int_0^R u''(t)dt\Big)ds + \int_R^r \Big(\int_R^s u''(t)dt\Big)ds \\ &\leq u(R) + (r-R)\int_0^R u''(t)dt + C_2\int_R^r \Big[\int_R^s \Big(t\int_0^t u^{-1/3}(\tau)d\tau\Big)dt\Big]ds \\ &\leq u(R) + r\int_0^R u''(t)dt + C_2 r^3\int_0^r u^{-1/3}(s)ds. \end{aligned}$$

Consequently,

$$u(r) \leq C_3 r^3 \int_0^r u^{-1/3}(s)ds$$

for some $C_3 > 0$ and for $r$ large enough. Denote

$$H(r) := \int_0^r u^{-1/3}(s)ds.$$

Then the previous inequality can be rewritten as

$$H'(r)H^{1/3}(r) \geq C_4 r^{-1}$$

for some $C_4 > 0$ and for $r$ large enough. Integrating by parts leads us to

$$H(r) \geq C_5(\log r)^{3/4}$$

for some $C_5 > 0$ and for $r$ large. We now again make use of the representation

$$u''(r) = u''(0) + \int_0^r H(s)ds$$

to deduce that

$$u''(r) \geq \int_{r/2}^r H(s)ds.$$

From this we obtain $u''(r) \geq (r/2)H(r/2)$ which helps us to conclude that

$$u''(r) \geq C_6 r(\log r)^{3/4}$$

for some $C_6 > 0$ provided $r$ is large enough. Upon repeating this trick one more time, we deduce that

$$u(r) = u(0) + \int_0^r \Big(\int_0^s u''(t)dt\Big)ds \geq C_7 r^3(\log r)^{3/4}$$

for some $C_7 > 0$ and for $r$ sufficiently large. Hence,

$$u^{(4)}(r) = u^{-1/3}(r) \leq C_7 r^{-1}(\log r)^{-1/4}$$

for $r$ large and, by integrating four times, we arrive at

$$u(r) \le C_8 r^3 (\log r)^{3/4}$$

for some $C_8 > 0$ and for $r$ large. Hence, we have just shown that

$$C_7 r^3 (\log r)^{3/4} \le u(r) \le C_8 r^3 (\log r)^{3/4}$$

for $r$ large. Keep in mind that

$$rH(r) \ge C_5 r(\log r)^{3/4}$$

and for some $C_9 > 0$ and that

$$\int_0^r s u^{-1/3}(s) ds \le C_9 r(\log r)^{-1/4}$$

for $r$ large. Therefore, we deduce from the representation

$$u''(r) = u''(0) + rH(r) - \int_0^r s u^{-1/3}(s) ds$$

that

$$\lim_{r \nearrow +\infty} \frac{u''(r)}{rH(r)} = 1.$$

Using the l'Hôpital rule and noting that $u^{-1/3}(r) = o(r^{-1})$, we have

$$\lim_{r \nearrow +\infty} \frac{u'(r)}{r^2 H(r)} = \lim_{r \nearrow +\infty} \frac{u''(r)}{2rH(r) + r^2 u^{-1/3}(r)} = \lim_{r \nearrow +\infty} \frac{u''(r)}{2rH(r)} = \frac{1}{2}.$$

Hence, we again apply the l'Hôpital rule to get

$$\lim_{r \nearrow +\infty} \frac{u(r)}{r^3 H(r)} = \lim_{r \nearrow +\infty} \frac{u'(r)}{3r^2 H(r) + r^3 u^{-1/3}(r)} = \lim_{r \nearrow +\infty} \frac{u'(r)}{3r^2 H(r)} = \frac{1}{6}. \tag{3.21}$$

Replacing $u(r) = (H'(r))^{-3}$, we deduce from (3.21) that

$$\lim_{r \nearrow +\infty} \frac{1}{rH'(r) H^{1/3}(r)} = \Big(\frac{1}{6}\Big)^{1/3}.$$

It follows from the l'Hôpital rule that

$$\lim_{r \nearrow +\infty} \frac{4 \log r}{3H^{4/3}(r)} = \Big(\frac{1}{6}\Big)^{1/3}. \tag{3.22}$$

Combining (3.21) with (3.22), we have

$$\lim_{r \nearrow +\infty} \frac{u(r)}{r^3 (\log r)^{3/4}} = \frac{1}{6^{3/4}} \Big(\frac{4}{3}\Big)^{3/4} = \Big(\frac{2}{9}\Big)^{3/4}$$

as claimed. □

Now we consider the remaining case $\alpha < -1/3$.

**Proposition 3.8.** *Assume that $n = 1$ and $\alpha < -1/3$. Let $u$ be a positive radial solution to the problem (1.1) in $\mathbf{R}^n$, then $u$ has the following asymptotic behavior*

$$u(r) \sim \frac{1}{6} \Big( \int_0^{+\infty} u^\alpha(t) dt \Big) r^3.$$

*Proof.* By the same argument as in the proof of (3.7), it is not hard to see that there exists some constant $C > 0$ such that

$$u(r) \geq C r^{4/(1-\alpha)}$$

for any $r$ large enough. From this and $\alpha < -1/3$, we deduce that

$$N := \int_0^{+\infty} u^\alpha(t)dt < +\infty.$$

Keep in mind that $u$ is an even function; hence $u'(0) = u^{(3)}(0) = 0$. Therefore, it follows from the equation satisfied by $u$ that

$$u^{(3)}(r) = \int_0^r u^\alpha(t)dt.$$

Consequently, there holds

$$\lim_{r \nearrow +\infty} u^{(3)}(r) = N.$$

We are now in position to apply the l'Hôpital rule to get

$$\lim_{r \nearrow +\infty} \frac{u(r)}{r^3} = \lim_{r \nearrow +\infty} \frac{u^{(3)}(r)}{6} = \frac{N}{6}.$$

Thus

$$u(r) \sim \frac{r^3}{6} \int_0^{+\infty} u^\alpha(t)dt$$

as claimed. □

## Appendix A. Asymptotic behavior at infinity of positive radial solutions to $\Delta u = u$ in $\mathbf{R}^n$

In this appendix, we establish the asymptotic behavior at infinity of positive radial solutions to $\Delta u = u$ in $\mathbf{R}^n$, which was not mentioned in the work [YG05]. To achieve the goal, we follow the proof of Proposition 3.1 and prove the following simple result which is more or less known to experts in the field.

**Proposition A.1.** *Let $n \geq 1$ and $u$ be a positive solution to the ODE*

$$u'' + \frac{n-1}{r} u' = u \quad \textit{for } r > 0$$

*with $u'(0) = 0$. Then $u$ has the following asymptotic behavior*

$$u(r) \sim u(0) \Gamma\Big(\frac{n}{2}\Big) 2^{\frac{n-3}{2}} \pi^{-1/2} r^{-\frac{n-1}{2}} e^r$$

*near infinity.*

*Proof.* **The case** $n = 1$. In this case, our equation simply becomes $u''(r) = u(r)$ and the general solution to the above ODE is of the form

$$u(r) = C_1 e^r + C_2 e^{-r},$$

where $C_1$ and $C_2$ are constants given by

$$C_1 = \frac{1}{2}\big[u(0) + u'(0)\big], \quad C_2 = \frac{1}{2}\big[u(0) - u'(0)\big].$$

From this we obtain the following

$$u \sim \frac{u(0)}{2} e^r.$$

**The case** $n \geq 2$. In this case the general solution to the ODE is of the form

$$u(r) = C_1 u_1(r) + C_2 u_2(r),$$

where $u_1$ and $u_2$ are already given in the proof of Proposition 3.1, namely

$$u_1(r) = r^{-\frac{n}{2}+1} I_{\frac{n}{2}-1}(r), \quad u_2(r) = r^{-\frac{n}{2}+1} K_{\frac{n}{2}-1}(r)$$

where $I$ and $K$ are the modified Bessel functions of the first and second kinds respectively. Recall the asymptotic behavior for the modified Bessel functions at infinity

$$I_{\frac{n}{2}-1}(r) \sim \frac{e^r}{\sqrt{2\pi r}}, \quad K_{\frac{n}{2}-1}(r) \sim \frac{\sqrt{\pi} e^{-r}}{\sqrt{2r}}.$$

From this we deduce that

$$u \sim C_1 u_1(r).$$

To estimate the constant $C_2$ we recall that

$$K_{\frac{n}{2}-1}(r) \sim_{r \searrow 0} \begin{cases} \log \dfrac{1}{r} & \text{if } n = 2, \\ \dfrac{1}{2}\Gamma\Big(\dfrac{n}{2} - 1\Big)\Big(\dfrac{r}{2}\Big)^{-\frac{n}{2}+1} & \text{if } n > 2, \end{cases}$$

Therefore, by working around zero and because $u(0)$ is finite, we deduce that $C_2 = 0$, namely

$$u(r) = C_1 u_1(r).$$

Hence $C_1 = u(0)/u_1(0)$. For the particular solution $u_1$, as in the proof of Proposition 3.1, we know that

$$u_1(0) = \frac{2^{-\frac{n}{2}+1}}{\Gamma(n/2)}, \quad u_1(r) \sim \frac{1}{\sqrt{2\pi}} r^{-\frac{n-1}{2}} e^r.$$

Hence

$$u(r) \sim u(0)\Gamma\Big(\frac{n}{2}\Big) 2^{\frac{n}{2}-1} \frac{1}{\sqrt{2\pi}} r^{-\frac{n-1}{2}} e^r.$$

This completes the proof. □

## Acknowledgments

QAN would like to acknowledge the support of the Tosio Kato Fellowship awarded in 2018 resulting in an earlier version [**NNP18**]. This work is supported by Vietnam National Foundation for Science and Technology Development (NAFOSTED) under grant number 101.02-2021.24 "On several real and complex differential operators" and grant number 101.02-2021.33 "Qualitative properties of partial differential equations and functional inequalities". Last but not least, we thank a referee for his/her careful reading and nice suggestion leading to the present version.

(Q.A. Ngô) University of Science, Vietnam National University, Hanoi, Vietnam, ORCID iD: 0000-0002-3550-9689

*Email address*: nqanh@vnu.edu.vn

(V.H. Nguyen) Department of Mathematics, FPT University, Hanoi, Vietnam

*Email address*: hoangnv47@fe.edu.vn

(Q. H. Phan) Institute of Research and Development, Duy Tan University, Da Nang, Vietnam

*Email address*: phanquochung@dtu.edu.vn